\documentclass[review]{elsarticle}
\usepackage{amsmath}
\usepackage{lineno,hyperref}
\usepackage{siunitx}
\usepackage{subfigure}
\usepackage[ruled,linesnumbered,vlined]{algorithm2e}
\usepackage{array}
\usepackage{tabularx}
\usepackage{threeparttable}
\usepackage{geometry}
\modulolinenumbers[1]

\biboptions{sort&compress}

\begin{document}
\begin{frontmatter}

\title{Efficient mesh deformation using radial basis functions with a grouping-circular-based greedy algorithm}

\author{Hong Fang}
\author{He Zhang}
\author{Fanli Shan\corref{mycorrespondingauthor}}
\cortext[mycorrespondingauthor]{Corresponding author at: Science and Technology on Space Physics Laboratory, Beijing 100076, P.R. China.}
\ead{tacyon@163.com}
\author{Ming Tie}
\author{Xing Zhang}
\author{Jinghua Sun}
\address{Science and Technology on Space Physics Laboratory, Beijing 100076, China.}

\begin{abstract}
A grouping-circular-based (GCB) greedy algorithm is proposed to promote the efficiency of mesh deformation. By incorporating the multigrid concept that the computational errors on the fine mesh can be approximated with those on the coarse mesh, this algorithm stochastically divides all boundary nodes into $m$ groups and uses the locally maximum radial basis functions (RBF) interpolation error of each group as an approximation to the globally maximum one of all boundary nodes in each iterative procedure for reducing the RBF support nodes. For this reason, it avoids the interpolation conducted at all boundary nodes and thus reduces the corresponding computational complexity from $O\left({N_c^2{N_b}} \right)$ to $O\left( {N_c^3} \right)$, where $N_b$ and $N_c$ denote the numbers of boundary nodes and support nodes, respectively. Besides, after $m$ iterations, the interpolation errors of all boundary nodes are computed once, thus allowing all boundary nodes can contribute to error control. Two canonical deformation problems of the ONERA M6 wing and the DLR-F6 Wing-Body-Nacelle-Pylon configuration are computed to validate the GCB greedy algorithm. The computational results show that the GCB greedy algorithm is able to remarkably promote the efficiency of computing the interpolation errors in the data reducing procedure by dozens of times. Because an increase of $m$ results in an increase of $N_c$, an appropriate range of $\left[ {{N_b}/{N_c},{\rm{ }}2{N_b}/{N_c}}\right]$ for $m$ is suggested to prevent too much additional computations for solving the linear algebraic system and computing the displacements of volume nodes induced by the increase of $N_c $. The results also show that the GCB greedy algorithm tends to generate a more significant efficiency improvement for mesh deformation when a larger-scale mesh is applied. Furthermore, this algorithm can produce a deformed mesh with a comparable quality to the undeformed one and retain the grid orthogonality and grid spacing near the solid surface for both structured and unstructured meshes.

\end{abstract}

\begin{keyword}
Grouping-circular\sep Greedy algorithm\sep Mesh deformation
\end{keyword}

\end{frontmatter}

\section{Introduction}\label{sec:1}
The essential issue for aerodynamic shape optimizations\cite{1, 2, 3,4}, aircraft icing simulations\cite{5,6} and aeroelasticity predictions\cite{7,8,9,10} is to allow the computational mesh deformed so as to facilitate reproductions of aerodynamic effects induced by shape shifting. To address this issue, various methods for mesh deformation have been developed in prior works, such as the spring analogy method\cite{11} and its improvements\cite{12,13}, elasticity analogy method\cite{14,15,16}, partially differential equation (PDE)-based method\cite{17,18}, optimization-based method\cite{19,20}, inverse distance weighted (IDW) method\cite{21}, Delaunay graph method\cite{22} and radial basis functions (RBF) method\cite{23,24,25,26} etc. Among these methods, the RBF one, first proposed by Boer et al.\cite{23}, has been widely adopted for mesh deformation in the past decade due to its simplicity, robustness and achievable mesh quality. In this method, the displacement of an arbitrary volume node in the calculation domain is determined by using the RBF interpolation, which is a weighted summation of a set of the RBFs. Each RBF is associated with a support node that is essentially a boundary node, and it can characterize the relative displacement between the volume node and the support node according to their radius vectors. This indicates that no connectivity information is required and the deformations of structured and unstructured meshes can be therefore treated in a uniform manner.

However, the use of the RBF method poses a challenge that it tends to generate an extensive data manipulation and computing, which results in a great amount of computational expenditure, especially when a large-scale mesh with a large number of boundary nodes is applied. As a result, an employment of data reduction for this method is required to overcome this challenge. Rendall and Allen\cite{24,25} proposed a greedy algorithm by minimizing the number of support nodes for data reduction. This algorithm was devised with an incorporation of an error-driven data reducing procedure, which starts from constructing the set of support nodes based on an arbitrary initial subset of boundary nodes, then solving the linear algebraic system to obtain the weighting coefficients and computing the RBF interpolation errors for all boundary nodes, after which adding the node with the maximum interpolation error to the set of support nodes. This procedure is performed repetitively until the maximum interpolation error is within allowance and an optimum reduced set of support nodes is then obtained. By using the greedy algorithm, computational expenditure can be remarkably decreased, thus leading to a significantly improved efficiency of mesh deformation. Based on the concept of the greedy algorithm, Wang and Mian\cite{26} modified the data reducing procedure by introducing a multi-level subspace scheme which treats the errors of boundary nodes in present interpolation as the objects of next interpolation. This modified algorithm can control the dimension of the linear algebraic system in a small scale and promote the efficiency of data reduction. Wei et al.\cite{27} developed a peak-selection algorithm by using the error-peak nodes to construct the set of support nodes, which results in a promoted efficiency of data reduction because multiple nodes can be added at a time

Data reduction can also be achieved by reducing the number of volume nodes at which the displacements need to be determined. For instance, Fang et al.\cite{28} proposed an algorithm to allow the interpolation is conducted on the Cartesian background mesh rather than on the computational mesh. Xie and Liu\cite{29} reduced the data by limiting the deforming region with the use of a restricted wall function. Kedward et al.\cite{KEDWARD2017732} presented a multiscale RBF interpolation, aimed at reducing the influence domain of RBF.

It is noteworthy that the computational complexity for solving the linear algebraic system to obtain the RBF weighting coefficients in the data reducing procedure of the greedy algorithm is $N_c^4/24$. This implies that $N_c$ cannot be too large, otherwise the computation will be too expensive. On the other hand, if the number of support nodes is insufficient, the requirement of the interpolation accuracy may not be satisfied. To resolve these conflicting issues, an algorithm based on the incremental Lower-Upper (LU) decomposition scheme was employed by Selim et al.\cite{30}, which can reduce the computational complexity for solving the linear algebraic system to $4N_c^3/3$, thus noticeably promoting the efficiency of data reduction and allowing a use of a sufficient number of support nodes to meet the requirement of the interpolation accuracy. For the same purpose, Fang et al.\cite{FANG2019183} proposed an algorithm based on the recurrence Cholesky (RC) decomposition scheme to acquire the recurrence solution of the lower triangular matrix. This algorithm is able to reduce the computational complexity to as small as $2N_c^3/3$ and further promote the efficiency of data reduction.

It is also noted that the computational complexity for computing the RBF interpolation errors at boundary nodes in the data reducing procedure of the greedy algorithm is $3N_c^2N_b/2$ with $N_b$ representing the number of boundary nodes. It is obvious that too large   may result in an unaffordable computation. To overcome this problem, Strofylas et al.\cite{31} constructed a multigrid agglomeration algorithm to decrease the number of boundary nodes used in the data reducing procedure and thus promote the efficiency of data reduction. However, the implementation of this algorithm is complicated, which may limit its application to mesh deformation.

In order to reduce the computational complexity for computing the RBF interpolation errors by decreasing the number of boundary nodes used in the data reducing procedure with an easy implementation, a concept based on a stochastic grouping to boundary nodes is proposed, and a grouping-circular-based (GCB) greedy algorithm which incorporates this concept is accordingly constructed in this work. In the GCB greedy algorithm, all boundary nodes are stochastically divided into $m$ groups and each of groups is a subset of all boundary nodes, thus allowing the interpolation conducted in each group rather than at all boundary nodes in each iterative procedure. For this reason, it is able to make the computational complexity for computing the RBF interpolation errors reduced by $m$ times, which helps to promote the efficiency of data reduction. Besides, this algorithm uses the node with the maximum interpolation error of each group to construct the set of support nodes.

In the following sections, the theory of the RBF interpolation is first given and the GCB greedy algorithm is then introduced in detail. Finally, it is applied in two deformation problems of ONERA M6 wing and DLR-F6 Wing-Body-Nacelle-Pylon configuration for algorithm validation.

\section{RBF interpolation}

The RBF interpolation is shown as below.

\begin{equation}
F\left( {\vec r} \right) = \sum\limits_{i = 1}^{{N_c}} {{w_i}\varphi \left( {\left\| {\vec r - {{\vec r}_i}} \right\|} \right)}\label{eq:1}
\end{equation}
Here, $\vec r$ and $\vec r_i$ are the radius vectors of an arbitrary volume node in the calculation domain and the support node of the $i$th RBF, respectively; $\varphi$ is the general form of the RBF and it is a function of $\left\| {\vec r - {{\vec r}_i}} \right\|$ which stands for the distance between $\vec r$ and $\vec r_i$; $w_i$ is the weighting coefficient of the $i$th RBF and $F$ represents the weighted interpolation with a set of RBFs. The Wendland's $C^2$ function\cite{34,WENDLAND1998} as given in Eq. (\ref{eq:2}), which is recognized applicable to the RBF-based mesh deformation, is employed in this work.

\begin{equation}
\varphi \left( \eta  \right) = \left\{ {\begin{array}{*{20}{c}}
   {{{\left( {1 - \eta } \right)}^4}\left( {4\eta  + 1} \right),} &{\rm{for}\ \eta  \le 1}  \\
   {{\rm{\ \ \ \ \ \ \ \ \ \ \ \ \ \ \ \ \ \ \ \ }}0,} &{\rm{for}\ \eta  > 1}  \\
\end{array}} \right.
\label{eq:2}
\end{equation}
Here, $\eta$ is the non-dimensional distance normalized by the RBF radius. In the case that all boundary nodes are used to construct the set of support nodes, the following equation must be satisfied to ensure that the RBF interpolation holds for all boundary nodes.

\begin{equation}
F({\vec r_j}) = \Delta {\vec r_j} = \sum\limits_{i = 1}^{{N_b}} {{w_i}\varphi \left( {\left\| {{{\vec r}_j} - {{\vec r}_i}} \right\|} \right)} ,{\rm{\ \ \ }}j = 1,{\rm{\ }}2,{\rm{\ }} \cdots ,{\rm{\ }}{N_b}
\label{eq:3}
\end{equation}
Here, $N_b$ is the numbers of all boundary nodes; $\Delta {\vec r_j}$ is the displacement vector of the $j$th boundary node. In three dimensions, Eq. (\ref{eq:3}) can be rewritten in the matrix forms as below.

\begin{equation}
\begin{array}{l}
 {\bf{\Phi }}{{\bf{W}}^x} = \Delta {\bf{X}} \\
 {\bf{\Phi }}{{\bf{W}}^y} = \Delta {\bf{Y}} \\
 {\bf{\Phi }}{{\bf{W}}^z} = \Delta {\bf{Z}} \\
 \end{array}
\label{eq:4}
\end{equation}
Here $\Delta {\bf{X}} = {\left\{ {\Delta {x_1}{\rm{,\ }} \cdots {\rm{, }}\Delta {x_j}{\rm{,\ }} \cdots {\rm{,\ }}\Delta {x_{{N_b}}}} \right\}^T}$ and $\Delta {\bf{Y}} = {\left\{ {\Delta {y_1}{\rm{,\ }} \cdots {\rm{, }}\Delta {y_j}{\rm{,\ }} \cdots {\rm{,\ }}\Delta {y_{{N_b}}}} \right\}^T}$ and $\Delta {\bf{Z}} = {\left\{ {\Delta {z_1}{\rm{,\ }} \cdots {\rm{, }}\Delta {z_j}{\rm{,\ }} \cdots {\rm{,\ }}\Delta {z_{{N_b}}}} \right\}^T}$ are the displacement components of boundary nodes in the $x$, $y$ and $z$ directions, respectively; $\mathbf{\Phi }$ is the RBF matrix and it can be expressed in the following form.

\begin{equation}
\begin{array}{l}
{\bf{\Phi }} = \left[ {\begin{array}{*{20}{c}}
   {\varphi \left( {\left\| {{{\vec r}_1} - {{\vec r}_1}} \right\|} \right)} &  \cdots  & {\varphi \left( {\left\| {{{\vec r}_1} - {{\vec r}_i}} \right\|} \right)} &  \cdots  & {\varphi \left( {\left\| {{{\vec r}_1} - {{\vec r}_{{N_b}}}} \right\|} \right)}  \\
    \vdots  &  \vdots  &  \vdots  &  \vdots  &  \vdots   \\
   {\varphi \left( {\left\| {{{\vec r}_j} - {{\vec r}_1}} \right\|} \right)} &  \cdots  & {\varphi \left( {\left\| {{{\vec r}_j} - {{\vec r}_i}} \right\|} \right)} &  \cdots  & {\varphi \left( {\left\| {{{\vec r}_j} - {{\vec r}_{{N_b}}}} \right\|} \right)}  \\
    \vdots  &  \vdots  &  \vdots  &  \vdots  &  \vdots   \\
   {\varphi \left( {\left\| {{{\vec r}_{{N_b}}} - {{\vec r}_1}} \right\|} \right)} &  \cdots  & {\varphi \left( {\left\| {{{\vec r}_{{N_b}}} - {{\vec r}_i}} \right\|} \right)} &  \cdots  & {\varphi \left( {\left\| {{{\vec r}_{{N_b}}} - {{\vec r}_{{N_b}}}} \right\|} \right)}  \\
\end{array}} \right]
 \end{array}
\label{eq:5}
\end{equation}
With the use of Eq. (\ref{eq:2}), ${\bf{\Phi }}$ becomes a symmetric and positive definite matrix\cite{34}, which indicates that the Cholesky decomposition scheme can be employed to solve the linear algebraic system as shown in Eq. (\ref{eq:4}), with a promoted efficiency to obtain the weighting coefficient vectors of ${{\bf{W}}^x} = $${\left\{ {w_1^x{\rm{,\ }} \cdots {\rm{,\ }}w_j^x{\rm{,\ }} \cdots {\rm{,\ }}w_{{N_b}}^x} \right\}^T}$ and ${{\bf{W}}^y} = $${\left\{ {w_1^y{\rm{,\ }} \cdots {\rm{,\ }}w_j^y{\rm{,\ }} \cdots {\rm{,\ }}w_{{N_b}}^y} \right\}^T}$ and ${{\bf{W}}^z} = $${\left\{ {w_1^z{\rm{,\ }} \cdots {\rm{,\ }}w_j^z{\rm{,\ }} \cdots {\rm{,\ }}w_{{N_b}}^z} \right\}^T}$. After these vectors are obtained, the displacements of any volume node in the calculation domain in three dimensions can be therefore determined by using the following equations, where ${N_v}{\rm{ }}$ denotes the number of volume nodes.

\begin{equation}
\begin{array}{l}
 \Delta {x_k} = \sum\limits_{i = 1}^{{N_b}} {w_i^x\varphi \left( {\left\| {{{\vec r}_k} - {{\vec r}_i}} \right\|} \right)}  \\
 \Delta {y_k} = \sum\limits_{i = 1}^{{N_b}} {w_i^y\varphi \left( {\left\| {{{\vec r}_k} - {{\vec r}_i}} \right\|} \right)} ,{\rm{\ \ \ }}k = 1,{\rm{\ }}2,{\rm{\ }} \cdots ,{\rm{\ }}{N_v}{\rm{\ \ }} \\
 \Delta {z_k} = \sum\limits_{i = 1}^{{N_b}} {w_i^z\varphi \left( {\left\| {{{\vec r}_k} - {{\vec r}_i}} \right\|} \right)}  \\
 \end{array}
\label{eq:6}
\end{equation}
It is clearly seen that the RBF interpolation based on Eq. (\ref{eq:6}) requires no connectivity information, thus avoiding the dependence on mesh topology. As a result, this method can be applied to mesh deformation with arbitrary mesh types, including structured mesh, unstructured mesh and hybrid mesh etc.

\section{GCB greedy algorithm}

Using all boundary nodes to construct the matrix ${\bf{\Phi }}$ is prohibitively expensive, especially for a three-dimensional large-scale mesh. For this reason, a greedy algorithm, proposed by Rendall and Allen\cite{24}, is used to minimize the dimension of ${\bf{\Phi }}$ by reducing the number of support nodes. However, as introduced in section \ref{sec:1}, the computational complexity for computing the RBF interpolation errors at boundary nodes in the data reducing procedure of this algorithm is $3N_c^2{N_b}/2$, which can result in a great amount of computation. Table 1 summarizes the computational complexities and the corresponding amounts of computation for different processes when performing mesh deformation with the use of the greedy algorithm for a typical mesh deformation problem which involves boundary layers. In this problem, the numbers of volume nodes, boundary nodes and support nodes are 4,000,000, 80,000 and 2,000, respectively. It is shown that the percentage of the total amount of computation resulted from computing the RBF interpolation errors is as high as 84.67\%, indicating that data reducing efforts should focus on this process to promote the efficiency of the greedy algorithm.

\begin{table}[ht!]\footnotesize
\caption{Computational complexities and amounts of computation when performing mesh deformation with the use of the greedy algorithm.}\label{tab.Compare}
\label{tab:1}
\begin{center}
\renewcommand{\arraystretch}{1.0}
\begin{tabularx}{16cm}{m{6.6cm}<{\centering}m{2.0cm}<{\centering}m{2.0cm}<{\centering}X<{\centering}}
\hline \noalign{\smallskip}
Process & Computational complexity & Amount of computation & Percentage of total amount of computation(\%) \\
\hline
Solving linear algebraic system with RC decomposition scheme 	& $2N_c^3/3$       & $5.3 \times {10^9}$  &0.93 \\
Computing RBF interpolation errors                          	&$3N_c^2N_b^2/2$   &$4.8 \times 10^{11}$  &84.67\\
Computin displacements of volume nodes                       	&   $10N_cN_v$     & $8.0 \times 10^{10}$ &14.11\\
Others                                                          &-                 &-                     &0.29 \\
\hline
\end{tabularx}
\medskip
\end{center}
\end{table}

In order to meet the requirement of the interpolation accuracy, the number of support nodes ($N_c$) should be sufficient. Therefore, the key to minimize the computational complexity is to decrease the number of boundary nodes used in the data reducing procedure ($N_b$). It is noted that the principle of the greedy algorithm is to construct a reduced set of support nodes by adding the node with the maximum interpolation error to the set in an iterative procedure of data reduction until the required interpolation accuracy is satisfied. For this reason, it is essentially an error-driven successive approximation algorithm and the exact computation for the maximum interpolation error is therefore not required. Furthermore, the sequence in which the support nodes are added does not affect the accuracy. Based on these considerations, the GCB greedy algorithm is constructed in this work. This algorithm treats the locally maximum interpolation error in a subset of all boundary nodes as an approximation to the globally maximum one to decrease $N_b$ and adds the node with the locally maximum interpolation error to the set of support nodes. It is similar to the multigrid concept which treats the computational errors on the coarse mesh as approximations to those on the fine mesh.

The procedure of performing the GCB greedy algorithm is given in Algorithm \ref{alg:1}.
Prior to mesh deformation, the set of all boundary nodes, denoted as $\textbf{\emph{G}}$, are stochastically divided into $m$ groups: $G_0$, $G_1$, $\cdots$, $G_{m-1}$. These groups satisfy the following conditions.
$$
\begin{array}{l}
\textbf{\emph{G}} = \left\{ {\left. {{G_i}} \right|i = 1,2, \cdots ,m - 1} \right\} \\
 {G_i} \cap {G_j} = \emptyset \begin{array}{*{20}{c}} {} & {i \ne j}  \\
\end{array}\begin{array}{*{20}{c}} {} & {} & {} &{i,j = 1,2, \cdots ,m - 1}  \\
\end{array} \\
 \left| {card({G_i}) - card({G_j})} \right| \le 1\begin{array}{*{20}{c}} {} & {i,j = 1,2, \cdots ,m - 1}  \\
\end{array} \\
\end{array}
$$
Here, $card$  represents the number of nodes of the group. In the $k$th iterative procedure of data reduction, only the $i$th ($i = k\bmod m$) group is active and the interpolation errors of the nodes in this group are computed, after which the node with the locally maximum interpolation error is added to the set of support nodes. After   iterative procedures are performed, the interpolation errors of all boundary nodes are computed once. As a result, all boundary nodes are used for error control, which facilitates ensuring the accuracy, thus leading to a superiority to other boundary node-based reducing algorithms which compute the errors only at selected boundary nodes. In addition, the stochastic grouping is helpful to approximate the globally maximum interpolation error of all boundary nodes with the locally maximum one of each group, which can accelerate the iteration. Another advantage of the GCB greedy algorithm is that it is easier to implement compared to other algorithms.

\begin{algorithm}[h]
\caption{The Grouping-Circular-Based(GCB) greedy algorithm}\label{alg:1}
\KwIn{$\textbf{\emph{G}}$,$\Delta {\bf{X}}$,$\Delta {\bf{Y}}$,$\Delta {\bf{Z}}$,$m$,$E^*$,$N_c^{max}$}
\KwOut{$\textbf{\emph{S}}$,$\textbf{W}^x$,$\textbf{W}^y$,$\textbf{W}^z$}
Divide $\textbf{\emph{G}}$ stochastically and uniformly into ${G_0},{G_1}, \cdots ,{G_i}, \cdots ,{G_{m - 1}}$\;
Select any 3 nodes $n_0,n_1,n_2$ from set $\textbf{\emph{G}}$ and add them to $\textbf{\emph{S}}$ (the set of support nodes)\;
\For{$k=3$ \emph{\KwTo} $N_c^{max}$ ${\bf step}$ $1$}{
    Compute ${\bf{\Phi }}$ for $\textbf{\emph{S}}$\;
    Get $\textbf{W}^x,\textbf{W}^y,\textbf{W}^z$ by solving ${\bf{\Phi}}(\textbf{W}^x,\textbf{W}^y,\textbf{W}^z)=(\Delta {\bf{X}},\Delta {\bf{Y}},\Delta {\bf{Z}})$\;
    $i$ = $k$ $\textbf {mod}$ $m$\;
    Compute the interpolation error in $G_i$\;
    Get the locally maximal interpolation error E and the corresponding node $n$ in $G_i$\;
    \If{$E > E^* $ }{
        App $n$ to \textbf{\emph{S}}
    }
}
\end{algorithm}

The GCB greedy algorithm is able to promote the efficiency of data reduction by reducing the computational complexity for computing the interpolation errors from $3N_c^2N_b/2$ of the traditional greedy algorithm to $3N_c^2N_b/{(2m)} $. It is noteworthy that this number will decrease to $2N_c^3/3 $ in the premise of $m=9N_b/(4N_c) $, which is equivalent to solve the linear algebraic system with the use of the RC decomposition scheme. Since the percentage of the total amount of computation for this process is less than 1\% as shown in Table \ref{tab:1}, it can be expected that the use of the GCB greedy algorithm is able to remarkably promote the efficiency of data reduction.

\section{Results and discussions}
In this section, two cases are employed to validate the GCB greedy algorithm. The first is the deformation of the ONERA M6 wing. The second is the deformation of the DLR-F6 Wing-Body-Nacelle-Pylon configuration.

\subsection{Deformation of ONERA M6 wing}\label{sec:4.1}
\begin{figure}[htp]
\begin{minipage}[t]{0.48\linewidth}
\centering
\includegraphics[height=6cm]{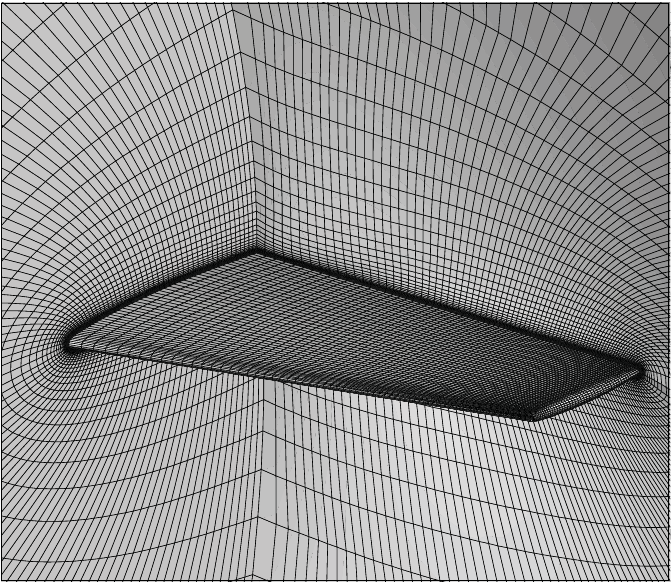}
\caption[1]{The schematic diagram of the ONERA M6 wing with its mesh distribution.}\label{fig:2}
\end{minipage}%
\hfill
\begin{minipage}[t]{0.48\linewidth}
\centering
\includegraphics[height=6cm]{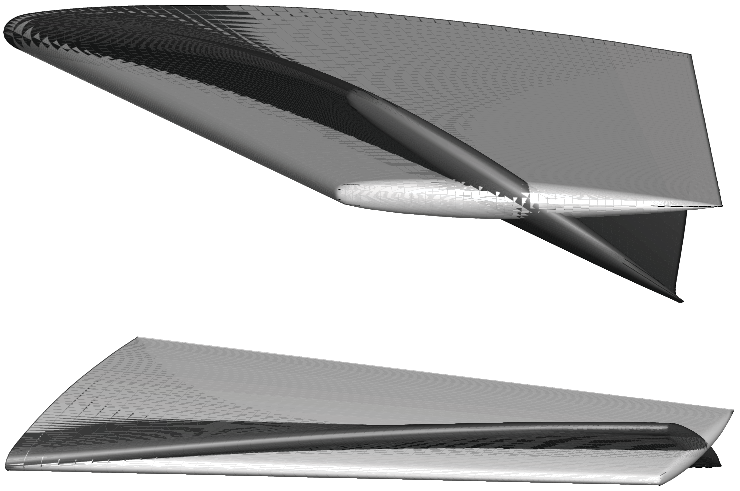}
\caption[1]{Shapes of the ONERA M6 wing before and after deformation.}\label{fig:3}
\end{minipage}%
\end{figure}

The schematic diagram of the ONERA M6 wing is illustrated in Fig. 1. A structured mesh, comprised of 1,310,720 elements and 1,353,105 nodes with 10,499 nodes on surface, is used to discretize the calculation domain and its distribution is also illustrated in this figure. To reproduce the boundary layer on surface, the height of the first layer mesh is set to \SI{1E-5}. A bending-twisting-coupled deformation, as described in Eqs. (\ref{eq:7}) and (\ref{eq:8}), is employed to validate the GCB greedy algorithm.

\begin{equation}
\Delta y = 0.05z\sin \left( {\frac{z}{{2b}}\pi } \right)
\label{eq:7}
\end{equation}

\begin{equation}
\theta  = {\theta _m}\sin \left( {\frac{z}{{2b}}\pi } \right)
\label{eq:8}
\end{equation}
Here, $b$ is the length of the root chord and equals to 0.805; ${\theta _m}$ is the wing twisting angle around the $1/4$ chord and is set to $\ang{30}$. The shapes of the wing before and after experiencing this coupled deformation are shown in Fig. \ref{fig:3}. The RBF radius is valued as 7 and the allowable error for the RBF interpolation is specified as $1 \times {10^{ - 6}}$.

\begin{figure}[htp]
\centering
\includegraphics[width=16cm]{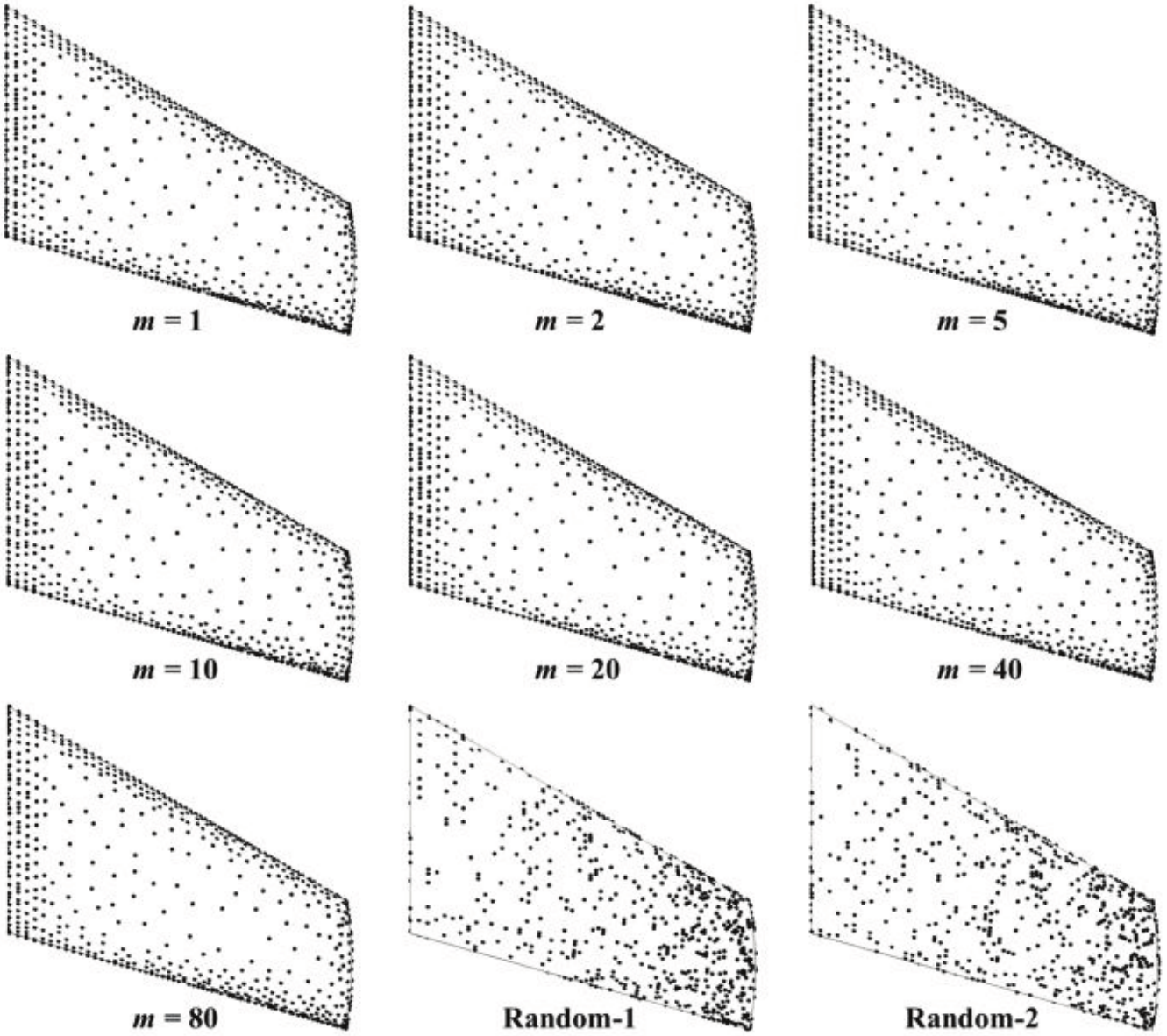}
\caption[1]{Distributions of support nodes for the deformation of the ONERA M6 wing ($N_c=1,104$)}\label{fig:baseNodeM6}
\end{figure}

The mesh deformation is first conducted using the traditional greedy algorithm and a total number of 1,104 nodes are selected as support nodes. For purpose of comparison, the number of support nodes is fixed at 1,104 when performing the GCB greedy algorithm. The distributions of support nodes on surface as $m$ varies from 1 to 80 are summarized in Fig. \ref{fig:baseNodeM6}. It is noted that the distributions of the GCB greedy algorithm are little deviated from that of the traditional greedy algorithm which is essentially a reduced GCB greedy algorithm in the case of $m=1$. On the contrary, the distributions of the two cases (Random-1 and Random-2) in which 1,104 support nodes are randomly selected are distinctly deviated from that of the traditional greedy algorithm. The deviations can be quantitatively evaluated by the Kullback-Leibler (KL) divergence which is defined in the following form.

\begin{equation}
KL\left( m_1||m_2 \right) = \sum \limits_{i=1}^{N_b}{d_i^{m_1} \log \frac{d_i^{m_1}}{d_i^{m_2}}} 
\label{eq:KLdef}
\end{equation}

\begin{equation}
d_i^{m_k} = \mathop {\min }\limits_{i \in G,{\rm{ }}j \in {S^{m_k}}} \left( {\left\| {{{\vec r}_i} - {{\vec r}_j}} \right\|} \right){\rm{,   }}k = 1,{\rm{ }}2{\rm{  }}
\end{equation}

Here, $S^{m_k}$ denotes the set of support nodes and $d_i^{m_k}$ describes the mapping relation between $S^{m_k}$ and $G$. According to Eq. \ref{eq:KLdef}, a smaller $KL$ divergence indicates the distributions of support nodes of two sets are more similar. Table \ref{tab:KL} lists the $KL$ divergences for different sets of supports nodes. It can be seen that the $KL$ divergences for the distributions of the GCB greedy algorithm are much smaller than those for the random cases, indicating that this algorithm is able to generate a more reasonable set of support nodes.

\begin{table}[ht!]\small
\caption{KL divergences for the deformation of the ONERA M6 wing ($N_c=1,1104$)}\label{tab:KL}
\begin{center}
\renewcommand{\arraystretch}{1.0}
\begin{tabularx}{8cm}{m{4cm}<{\centering}X<{\centering}}
\hline \noalign{\smallskip}
$S^{m_2}$ & $KL\left( {{m_1}||{m_2}} \right) \left( m_1=1 \right)$ \\
\hline
${m_2=2}$ 	      &182.5  \\
${m_2=5}$ 	      &183.8  \\
${m_2=10}$ 	      &182.1  \\
${m_2=20}$ 	      &183.6  \\
${m_2=40}$        &194.4  \\
${m_2=80}$        &189.3  \\
${Random-1}$      &1445.2 \\
${Random-2}$      &1404.3 \\
\hline
\end{tabularx}
\medskip
\end{center}
\end{table}

\begin{figure}[htp]
\centering
\includegraphics[width=9cm]{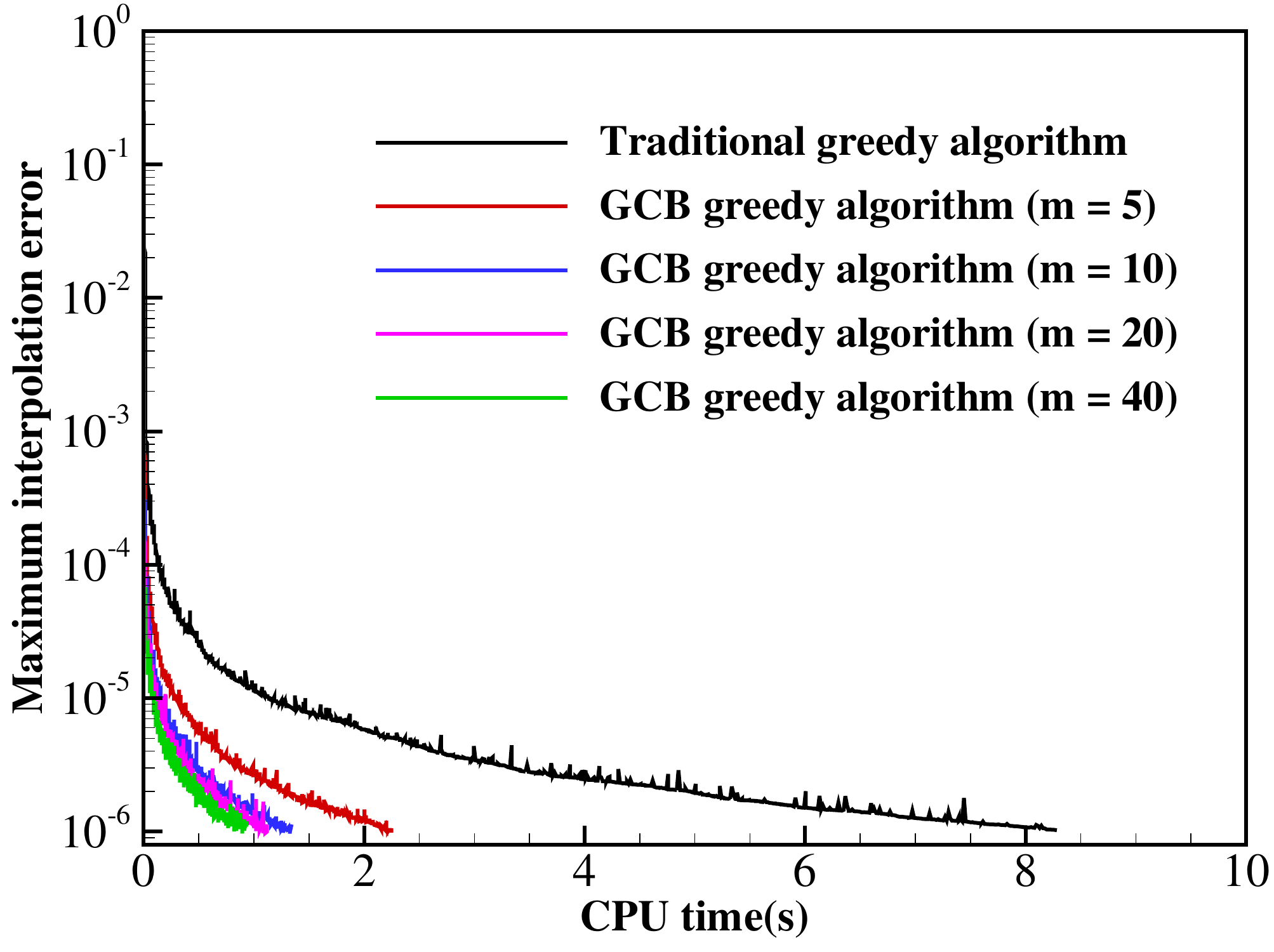}
\caption[1]{ Maximum interpolation error histories for the deformation of the ONERA M6 wing ($N_c$ = 1,104).}\label{fig:4}
\end{figure}

Fig. \ref{fig:4} displays the maximum interpolation error histories. For the traditional greedy algorithm, this error is globally determined among all boundary nodes, while for the GCB greedy algorithm, it is locally determined in the active group. It can be seen that the use of the GCB greedy algorithm results in a significantly reduced CPU time compared to that of the traditional greedy algorithm. Table \ref{tab:2} further manifests the superiority of the GCB greedy algorithm by listing the time consumptions for computing the interpolation errors ($t_1$) in the data reducing procedure. It is shown that $t_1$ decreases as the grouping number ($m$) increases. Compared to the traditional greedy algorithm, the GCB greedy algorithm with $m=40$ reduces $t_1$ from \SI{7.74}{s} to \SI{0.37}{s}, and the efficiency is therefore promoted by as high as 20.9 times. Considering that computing the interpolation errors using the traditional greedy algorithm takes over 80\% of the total amount of computation as shown in Table \ref{tab:1}, it is implied that the GCB greedy algorithm is able to overcome the challenge of the extensive data manipulation and computing posed by the traditional greedy algorithm. Table \ref{tab:2} also lists the time consumptions for solving the linear algebraic system ($t_1$). It is noted that $t_1$ and $t_2$ are of the same order of magnitude for $m \ge 10$, indicating that the corresponding computational complexity for computing the interpolation errors can decrease to $O\left( {N_c^3} \right)$ as long as $m$ satisfies certain condition.

\begin{table}[ht!]\small
\caption{Time consumptions for the deformation of the ONERA M6 wing ($N_c = 1,104$).}\label{tab:2}
\begin{center}
\renewcommand{\arraystretch}{1.0}
\begin{tabular}{cccc}
\hline \noalign{\smallskip}
Algorithm & $t_1(s)$ & $t_2(s)$ & $t_1+t_2(s)$ \\
\hline
Traditional greedy algorithm 	&7.74   &0.55   &8.29 \\
GCB greedy algorithm ($m = 5$) 	&1.69   &0.56   &2.25 \\
GCB greedy algorithm ($m = 10$) 	&0.75   &0.59   &1.34 \\
GCB greedy algorithm ($m = 20$)   &0.55   &0.57   &1.12 \\
GCB greedy algorithm ($m = 40$)   &0.37   &0.56   &0.93 \\
\hline
\end{tabular}
\medskip
\end{center}
\end{table}

\begin{table}[ht!]\small
\caption{Number of support nodes and time consumptions for the deformation of the ONERA M6 wing.}\label{tab:3}
\begin{center}
\renewcommand{\arraystretch}{1.0}
\begin{tabular}{cccccc}
\hline \noalign{\smallskip}
Algorithm &$N_c$  & $t_1(s)$ & $t_2(s)$  & $t_3(s)$ & $t_1+t_2+t_3(s)$ \\
\hline
Traditional greedy algorithm  &1104	&7.74   &0.55   &4.45   &12.74 \\
GCB greedy algorithm ($m = 5$)  &1107	&1.69   &0.57   &4.44   &6.70 \\
GCB greedy algorithm ($m = 10$) &1107	&0.76   &0.59   &4.49   &5.84 \\
GCB greedy algorithm ($m = 20$) &1137 &0.57   &0.58   &4.60   &5.75 \\
GCB greedy algorithm ($m = 40$) &1155 &0.41   &0.61   &4.67   &5.69 \\
GCB greedy algorithm ($m = 80$) &1202 &0.23   &0.67   &4.83   &5.73 \\
\hline
\end{tabular}
\medskip
\end{center}
\end{table}

\begin{figure}[htp]
\centering
\includegraphics[width=16cm]{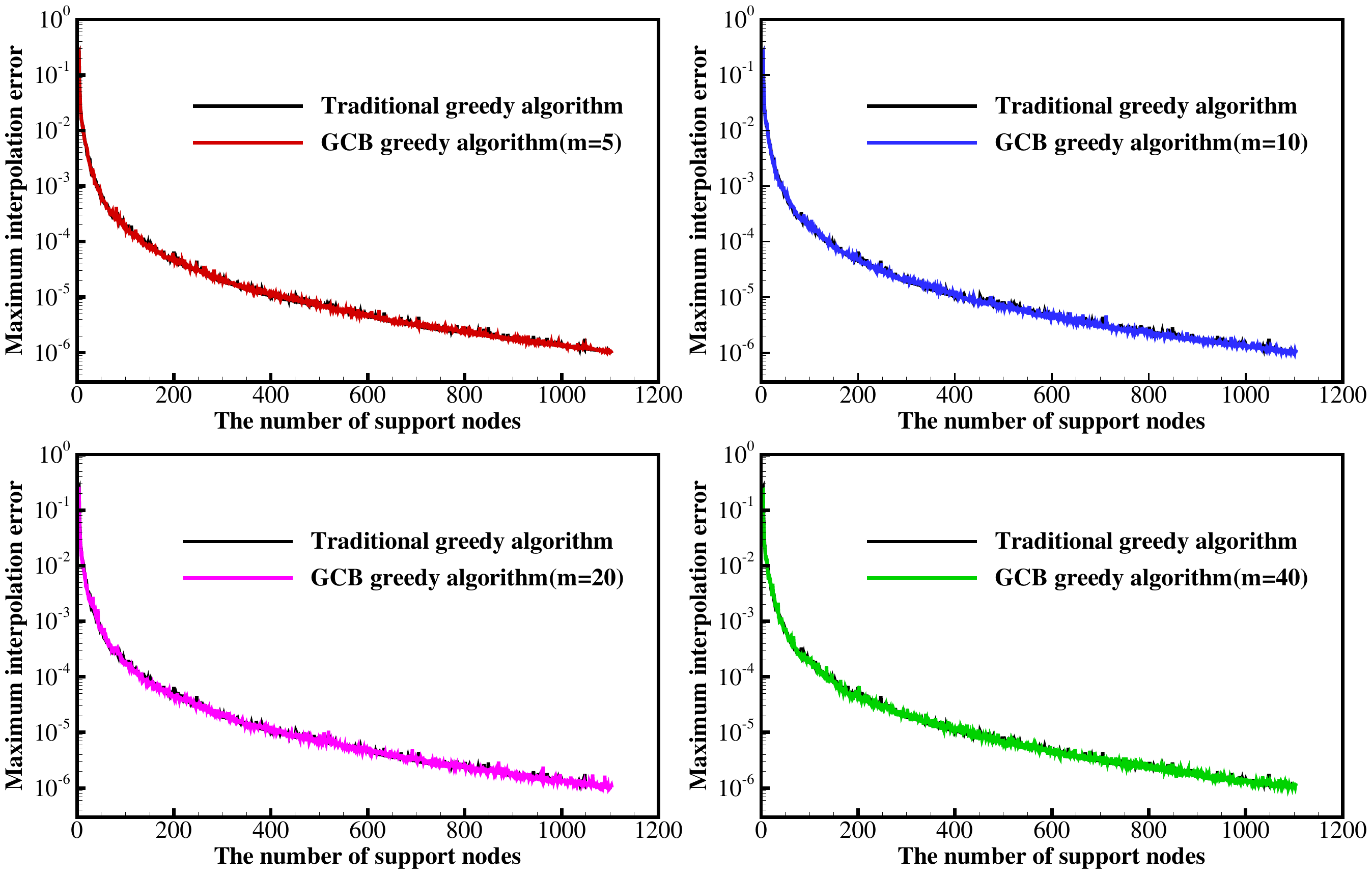}
\caption[1]{Maximum interpolation error histories in terms of the number of support nodes for the deformation of the ONERA M6 wing ($N_c$ = 1,104).}\label{fig:5}
\end{figure}

\begin{figure}[htp]
\centering
\begin{minipage}[t]{0.54\linewidth}
\centering
\includegraphics[height=5cm]{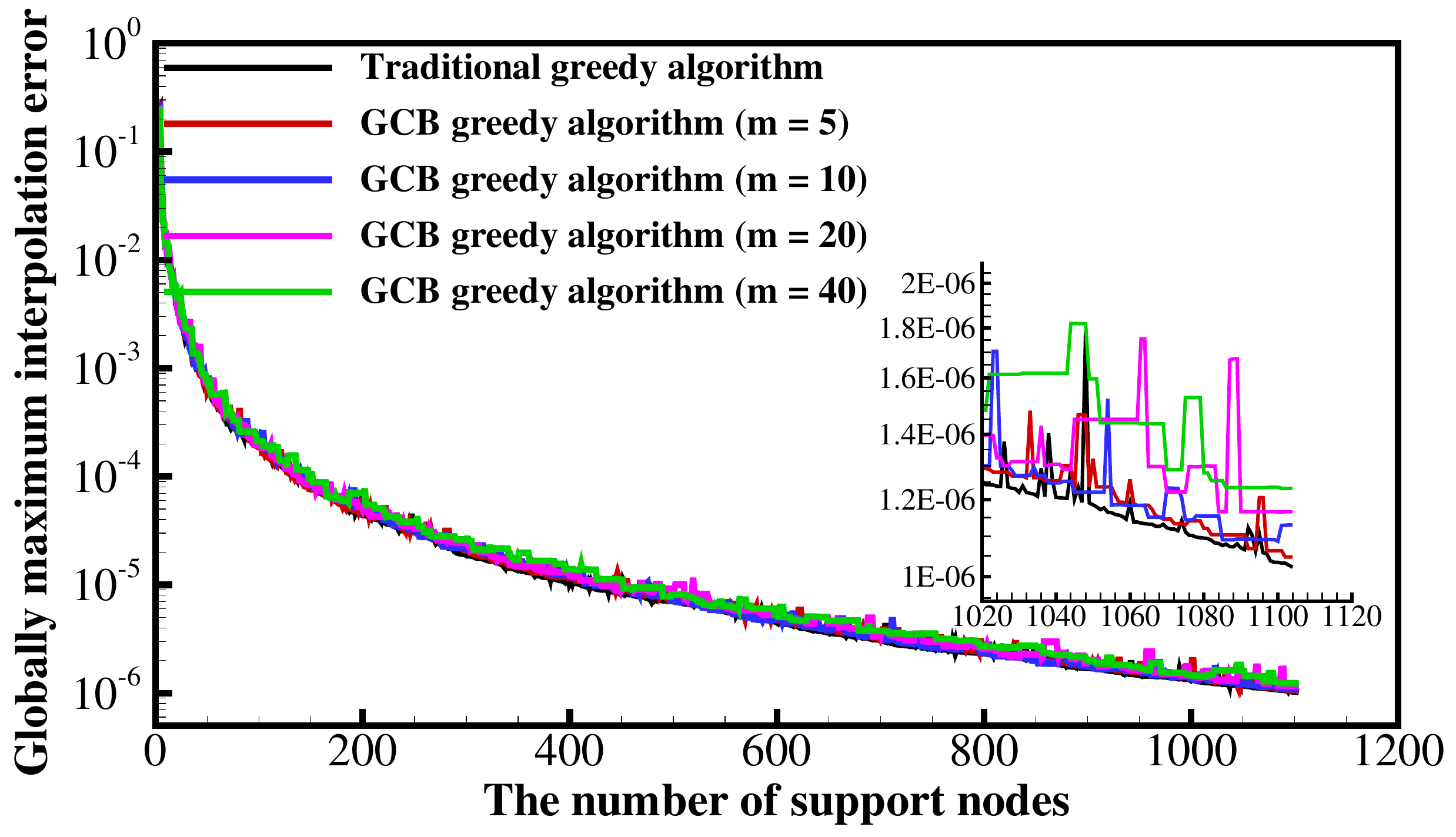}
\caption[1]{Globally maximum interpolation error histories in terms of the number of support nodes for the deformation of the ONERA M6 wing ($N_c$ = 1,104).}\label{fig:6}
\end{minipage}%
\hfill
\begin{minipage}[t]{0.42\linewidth}
\centering
\includegraphics[height=5cm]{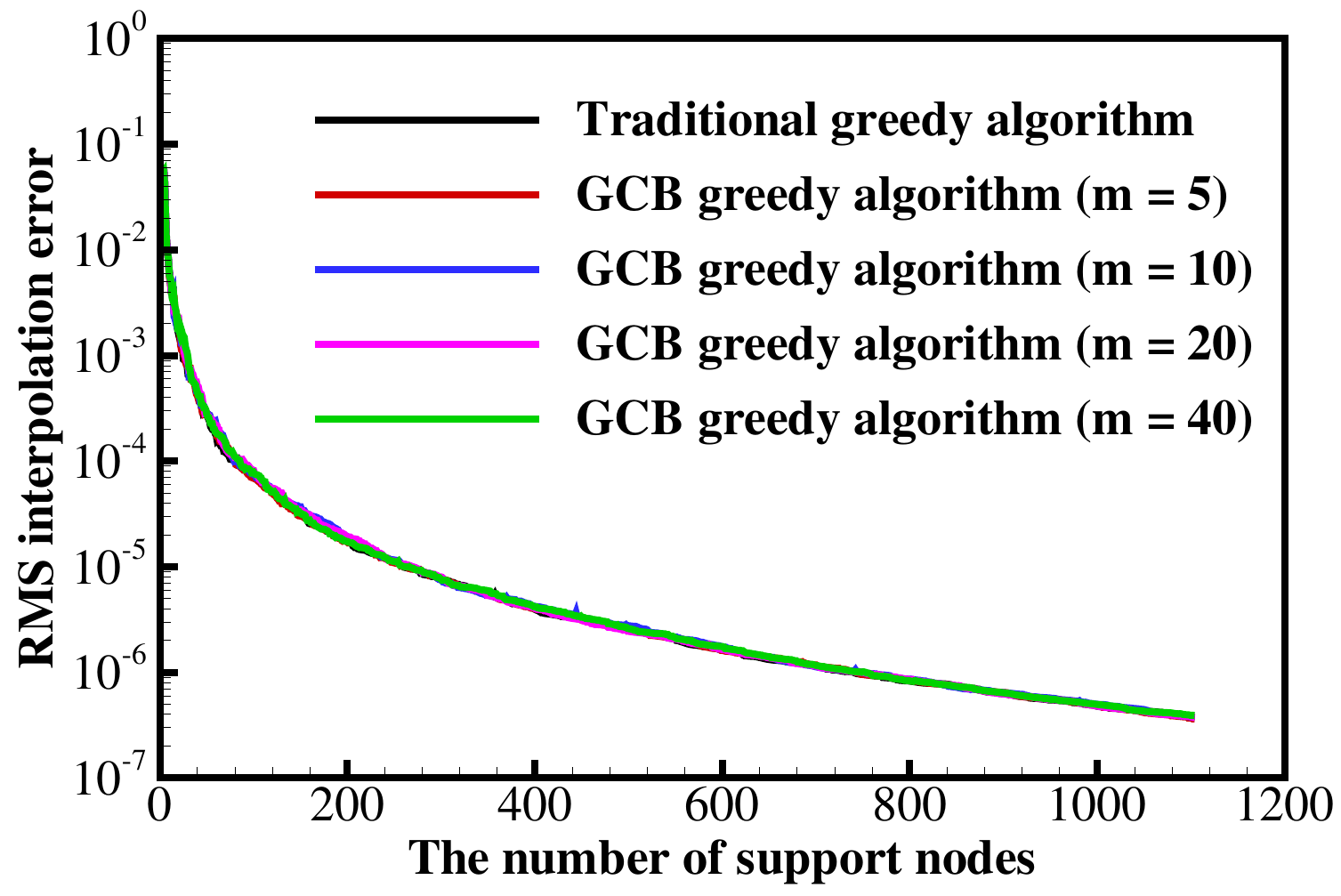}
\caption[1]{RMS interpolation error histories in terms of the number of support nodes for the deformation of the ONERA M6 wing ($N_c$ = 1,104).}\label{fig:7}
\end{minipage}%
\end{figure}

Given in \ref{fig:5} are the maximum interpolation error histories in terms of the number of support nodes. It is clearly seen that the errors of the GCB greedy algorithm decrease in a similar trend to that of the traditional greedy algorithm. This is associated with the similarity in the distribution of supports nodes as illustrated in Fig. \ref{fig:baseNodeM6}. Fig. \ref{fig:6} compares the globally maximum interpolation error histories of the two algorithms. It can be seen that the error decreasing trends of the two algorithms, on the whole, are in the similar manner, though slight deviations exist, especially in the case of $m \ge 20$. Fig. \ref{fig:7} further gives the root-mean-square (RMS) interpolation error histories for an overall comparison. It shows that the curves of the GCB greedy algorithm with $m$ varying from 5 to 40 are almost coincided with that of the traditional greedy algorithm. This indicates that the interpolation errors of each group are good approximations to those of all boundary nodes, thus validating the accuracy of the GCB greedy algorithm. After $m$ iterations, the interpolation errors of all boundary nodes are computed once, thus allowing all boundary nodes can contribute to error control. Other boundary node-based reducing algorithms, such as the multigrid agglomeration algorithm, however, compute the interpolation errors only at selected boundary nodes, which implies that the errors at unselected ones cannot be controlled.

The influences of $m$ on $N_c$ are also studied by removing the above fix on $N_c$ and is presented in Table \ref{tab:3}. It is noted that as $m$ grows to 80, $N_c$ is slowly increased by only 8.9\% from 1,104 to 1,202 in the premise of achieving the same interpolation accuracy (\SI{1E-6}). This can decrease the time consumption for computing the interpolation errors ($t_1$), but at a cost that it tends to increase the time consumptions for solving the linear algebraic system ($t_2$) and computing the displacements of volume nodes ($t_3$), which may reduce the efficiency of mesh deformation. This is because that the corresponding computational complexities for solving the linear algebraic system and computing the displacements of volume nodes increase as cubic and linear functions of $N_c$, respectively, as shown in Table \ref{tab:1}. For this reason, $m$ cannot be too large. It is generally appropriate to set $m$ in the range of $\left[ {{N_b}/{N_c},{\rm{\ }}2{N_b}/{N_c}} \right]$, which can ensure that the order of magnitude of the computational complexity for computing the interpolation errors is $O\left( {N_c^3} \right)$, and can also prevent too much additional computations induced by the increase of $N_c$. For present case, it is appropriate to set $m$ to 40.

\begin{figure}[htp]
\centering
\begin{minipage}[t]{0.54\linewidth}
\centering
\includegraphics[height=6cm]{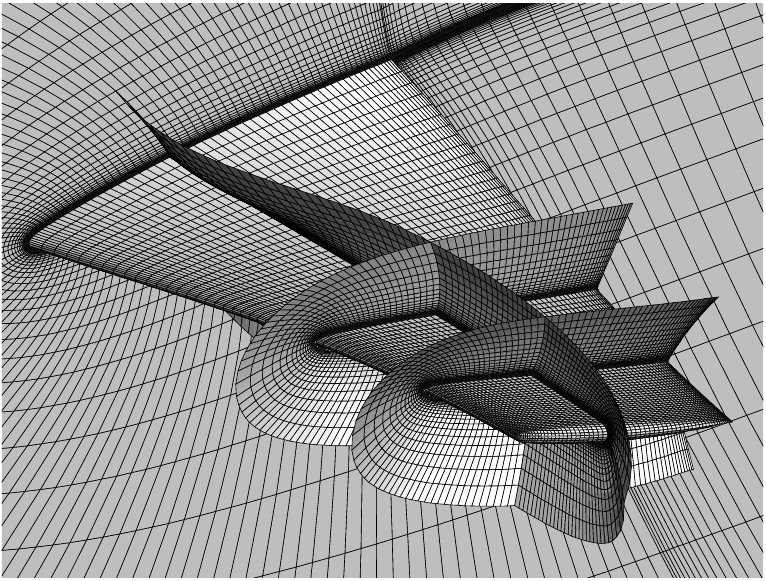}
\caption[1]{Mesh distribution of the ONERA M6 wing after deformation.}\label{fig:8}
\end{minipage}%
\hfill
\begin{minipage}[t]{0.42\linewidth}
\centering
\includegraphics[height=6cm]{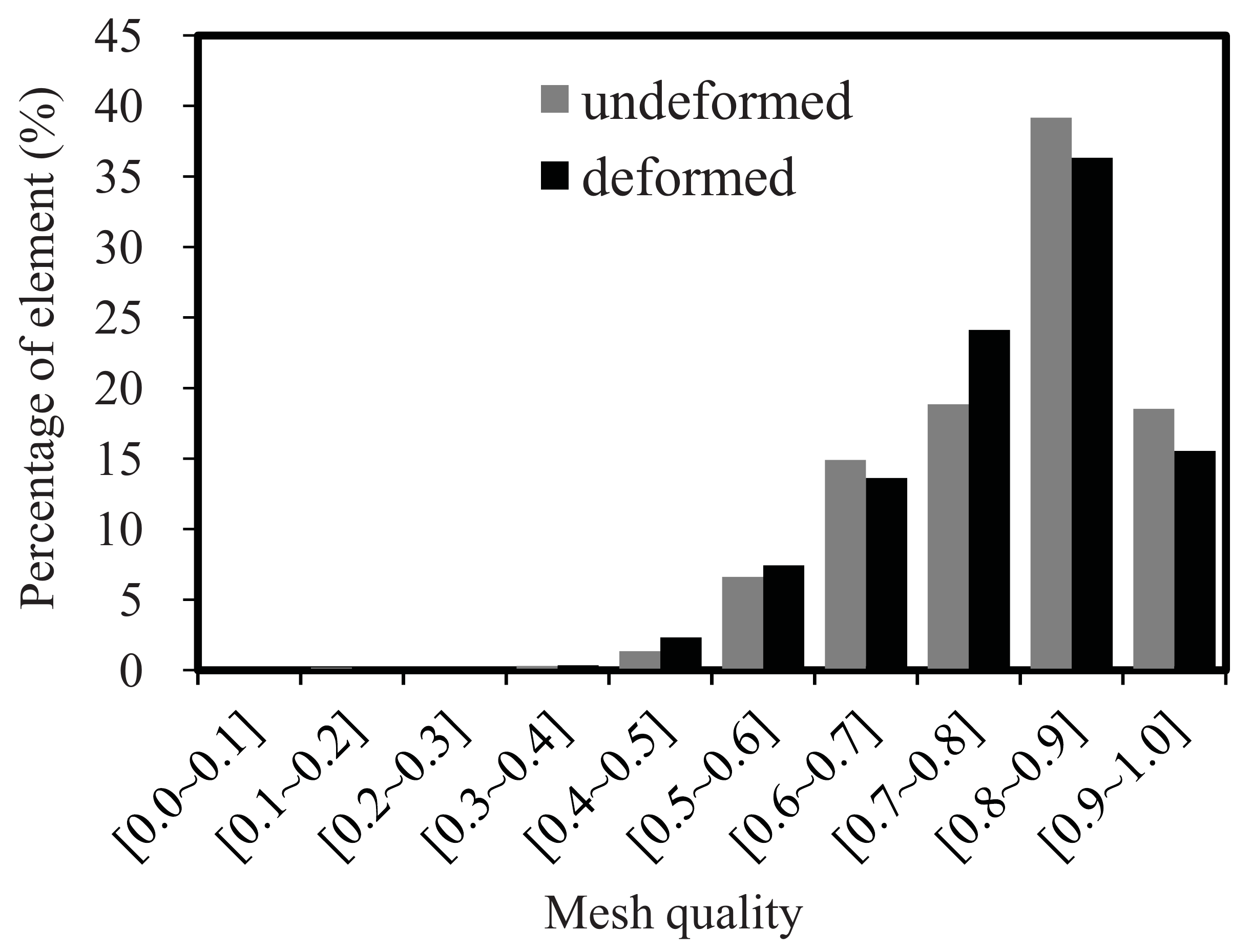}
\caption[1]{ Mesh qualities of the ONERA M6 wing undeformed and deformed.}\label{fig:mq1}
\end{minipage}%
\end{figure}

Mesh quality is another indicator to examine the algorithm. The mesh distribution after deformation for the GCB greedy algorithm is illustrated in Fig. \ref{fig:8} and the mesh quality defined in Eq. (\ref{eq:QualiDef}) \cite{NIU2017122} after deformation for this algorithm is shown in Fig. \ref{fig:mq1}.

\begin{equation}
q = 1 - \max \left( {\frac{{{\alpha _{\max }} - \alpha }}{{{{180}^ \circ } - \alpha }},\frac{{\alpha  - {\alpha _{\min }}}}{\alpha }} \right)
\label{eq:QualiDef}
\end{equation}
Here, $\alpha$ is the interior angle of the grid; ${\alpha _{\max}}$ and ${\alpha _{\min}}$ are the largest and smallest interior angles of the grid, respectively. If the value of $q$ tends to 1, the shape of the grid tends to a regular polygon and the mesh quality is considered fine. It is shown that the mesh qualities before and after deformation are comparable. In addition, the grid orthogonality and the grid spacing near surface are retained after deformation.

\subsection{Deformation of DLR-F6 Wing-Body-Nacelle-Pylon configuration}

\begin{figure}[htp]
\centering
\includegraphics[width=12cm]{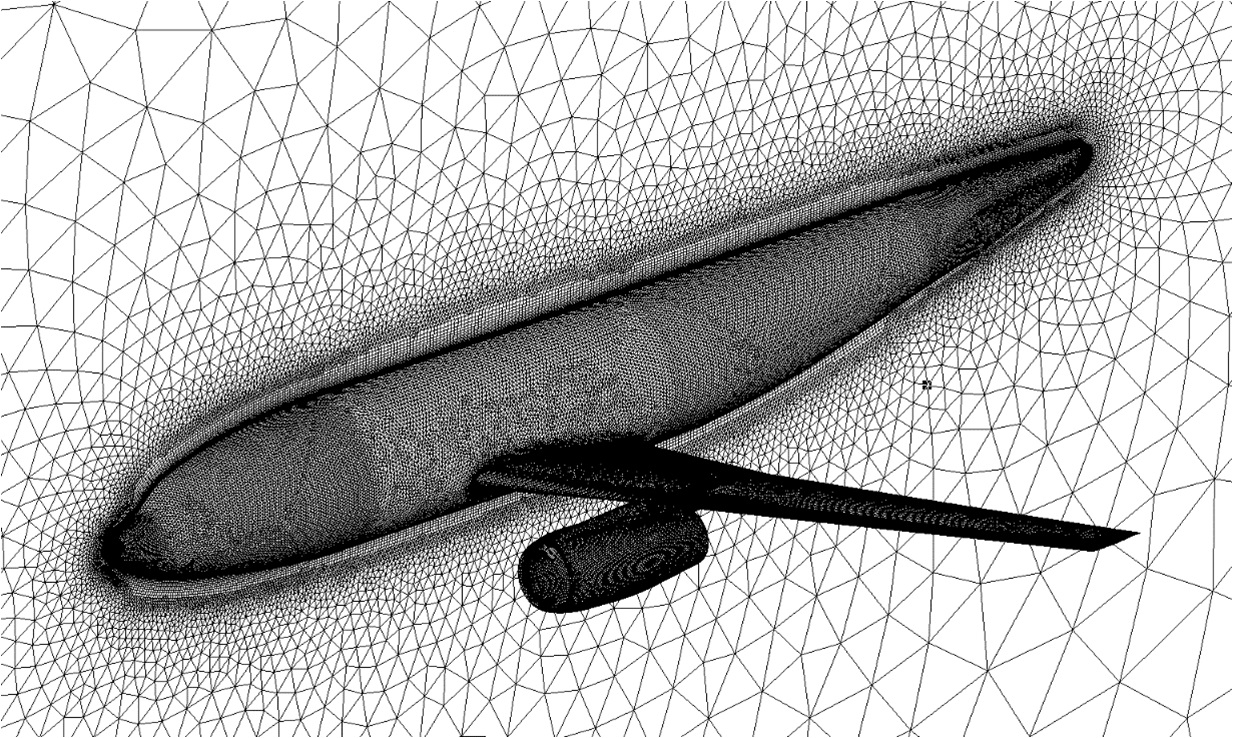}
\caption[1]{The schematic diagram of the DLR-F6 Wing-Body-Nacelle-Pylon configuration with its mesh distribution.}\label{fig:9}
\end{figure}

\begin{figure}[htp]
\centering
\includegraphics[width=12cm]{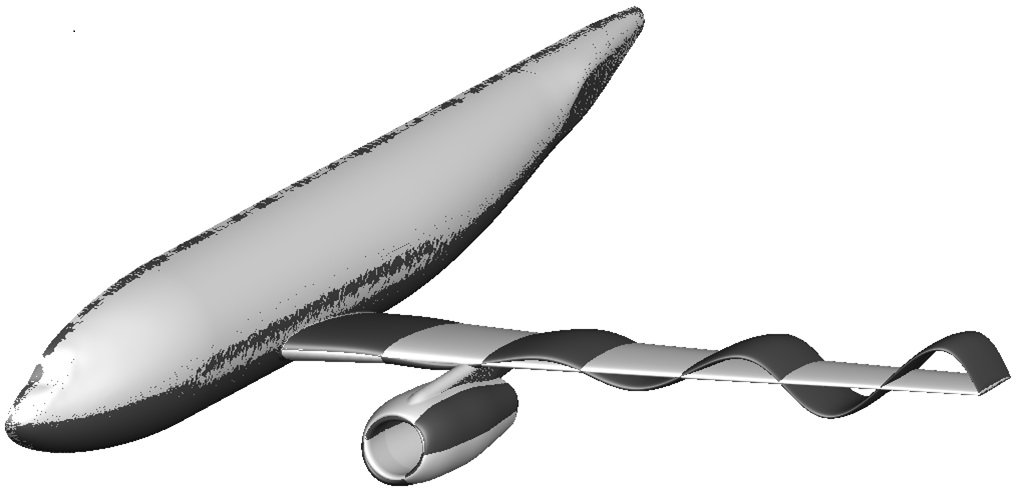}
\caption[1]{The DLR-F6 Wing-Body-Nacelle-Pylon configurations before and after deformation.}\label{fig:10}
\end{figure}

For present case, an unstructured mesh is applied to discretize the calculation domain. This mesh, as illustrated in Fig. \ref{fig:9}, is comprised of 9,233,536 elements and 4,007,571 nodes with 79,079 nodes on surface, which is larger than that of the case in section \ref{sec:4.1}. The height of the first layer mesh is also set to \SI{1E-5}. The configuration is deformed in a mode as below.

\begin{equation}
\Delta y = 0.3c\frac{{{z^2}}}{{{b^2}}}\sin \left( {8\pi \frac{z}{b}} \right)
\label{eq:10}
\end{equation}
Here, $b$ and $c$ are the span length and the mean aerodynamic chord length, respectively. A comparison of the configurations before and after deformation is given in Fig. \ref{fig:10}. The values of the RBF radius and the allowable error for the RBF interpolation are 5 and \SI{2E-6}, respectively.

\begin{figure}[htp]
\centering
\includegraphics[width=12cm]{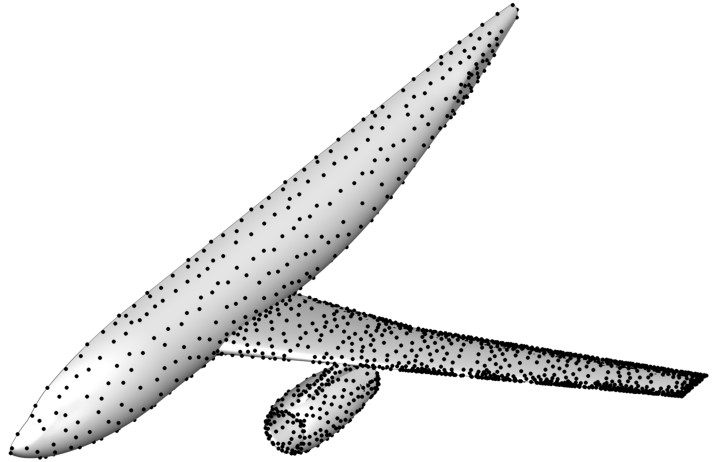}
\caption[1]{The distribution of support nodes for the deformation of the DLR-F6 Wing-Body-Nacelle-Pylon configuration ($N_c=2,520$)}\label{fig:baseNodeF6}
\end{figure}

\begin{figure}[htp]
\centering
\includegraphics[width=9cm]{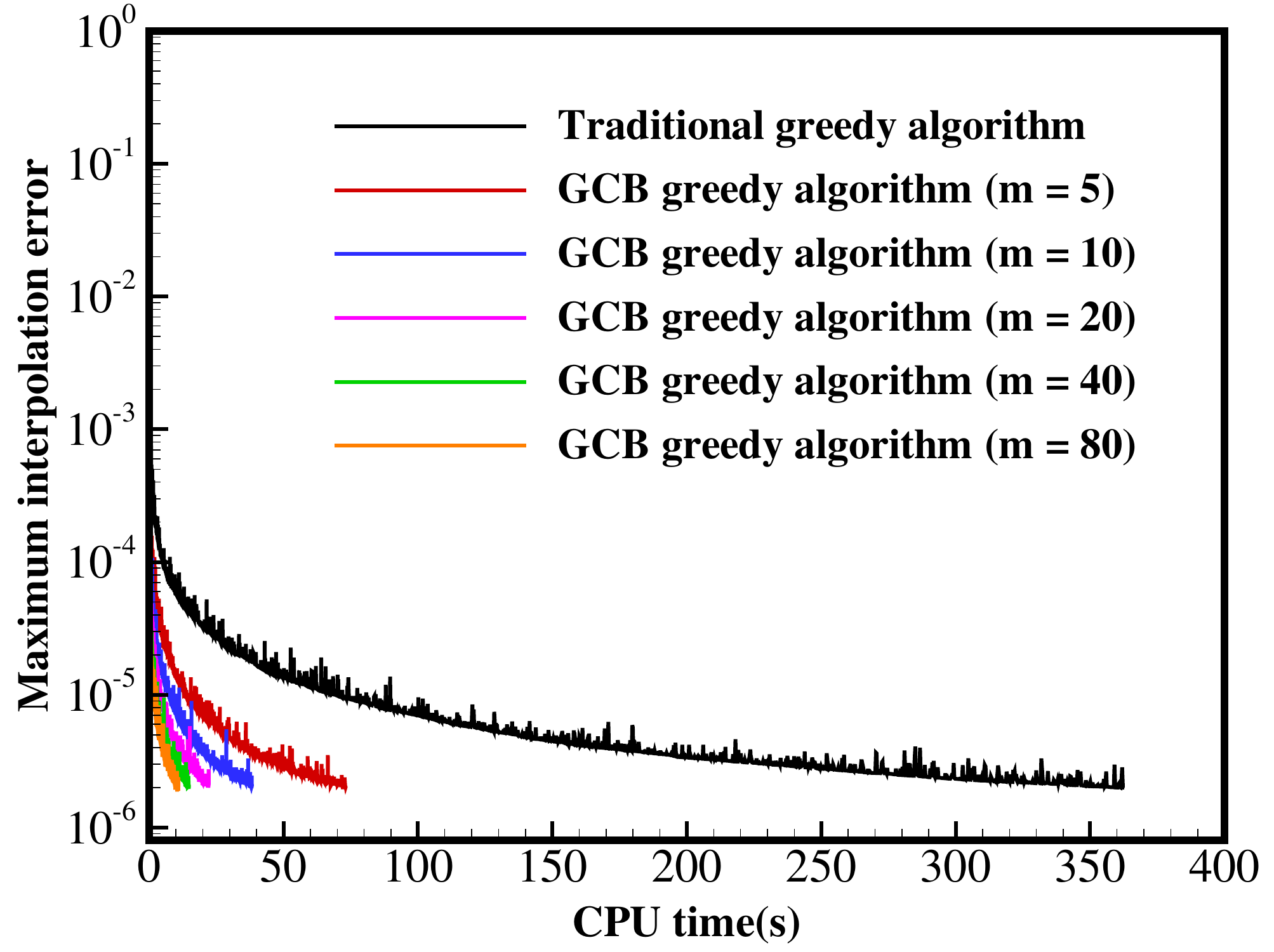}
\caption[1]{Maximum interpolation error histories for the deformation of the DLR-F6 Wing-Body-Nacelle-Pylon configuration ($N_c = 2,520$).}\label{fig:11}
\end{figure}

The number of support nodes selected by the traditional greedy algorithm is 2,520. This number is fixed when performing the GCB greedy algorithm for purpose of comparison. The distribution of support nodes generated by the traditional greedy algorithm is presented in Fig. \ref{fig:baseNodeF6}. Consistent with the case in section \ref{sec:4.1}, the GCB greedy algorithm can also generate similar distributions of support nodes to the traditional greedy algorithm for present case. Fig. \ref{fig:11} shows the maximum interpolation error histories of the traditional greedy algorithm and the GCB greedy algorithm with $m$ varying from 5 to 80. It can be seen that data reduction takes much more CPU time for a larger-scale mesh compared to the case in section \ref{sec:4.1} and the use of the GCB greedy algorithm is able to greatly reduce the CPU time. Table \ref{tab:5} quantitatively shows the magnitude of the decrease of the time consumptions in the data reducing procedure. For the traditional greedy algorithm ($m=1$), the sum of the time consumptions for computing the interpolation errors ($t_1$) and solving the linear algebraic system ($t_2$) are \SI{362.75}{s}, while for the GCB greedy algorithm with $m=80$, it decreases to \SI{10.95}{s}, indicating that the efficiency of data reduction is remarkably promoted by as high as 33.1 times. Since $t_2$ of the two algorithms are comparable, the efficiency promotion is mainly resulted from the employment of the grouping-circular concept.

\begin{table}[ht!]\small
\caption{Time consumptions for the deformation of the DLR-F6 Wing-Body-Nacelle-Pylon configuration ($N_c = 2,520$.)}\label{tab:5}
\begin{center}
\renewcommand{\arraystretch}{1.0}
\begin{tabular}{cccc}
\hline \noalign{\smallskip}
Algorithm & $t_1(s)$ & $t_2(s)$ & $t_1+t_2(s)$ \\
\hline
Traditional greedy algorithm 	&357.24   &5.51   &362.75 \\
GCB greedy algorithm ($m = 5$) 	&67.70    &5.49   &73.19 \\
GCB greedy algorithm ($m = 10$) 	&32.68    &5.58   &38.26 \\
GCB greedy algorithm ($m = 20$)   &16.86    &5.39   &22.25 \\
GCB greedy algorithm ($m = 40$)   &9.41     &5.44   &14.85 \\
GCB greedy algorithm ($m = 80$)   &5.58     &5.37   &10.95 \\
\hline
\end{tabular}
\medskip
\end{center}
\end{table}

\begin{figure}[htp]
\centering
\includegraphics[width=16cm]{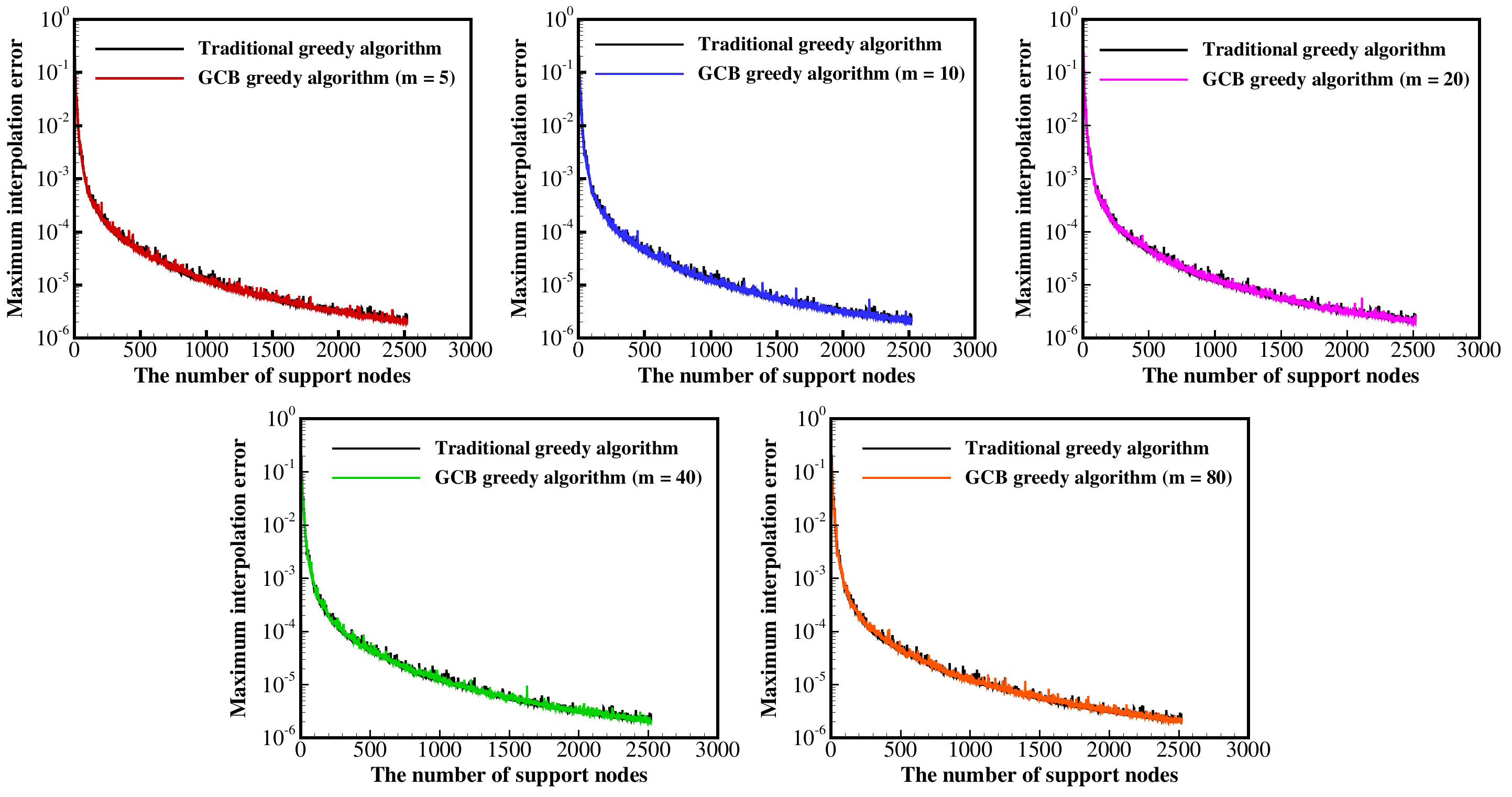}
\caption[1]{Maximum interpolation error histories in terms of the number of support nodes for the deformation of the DLR-F6 Wing-Body-Nacelle-Pylon configuration ($N_c = 2,520$).}\label{fig:12}
\end{figure}

The accuracy of the GCB greedy algorithm for present case can also be validated by examining the error histories in terms of the number of support nodes in comparison with the traditional greedy algorithm as summarized in Figs. \ref{fig:12} $\sim$ \ref{fig:15}. It is noted that the errors of the two algorithms generally follow similar decreasing trends, though as $m$ grows, the globally maximum interpolation error of the GCB greedy algorithm gradually deviates from that of the traditional greedy algorithm, and the order of magnitude of the deviations is $O\left( {{{10}^{ - 6}}} \right)$ as illustrated in Fig. \ref{fig:14}. It is also noted in Fig. \ref{fig:15} that the RMS interpolation error histories of the GCB greedy algorithm with $m$ varying from 5 to 80 are all in good agreement with that of the traditional greedy algorithm, demonstrating that using the interpolation errors of each group to approximate those of all boundary nodes is able to ensure the accuracy.

\begin{figure}[htp]
\centering
\includegraphics[width=16cm]{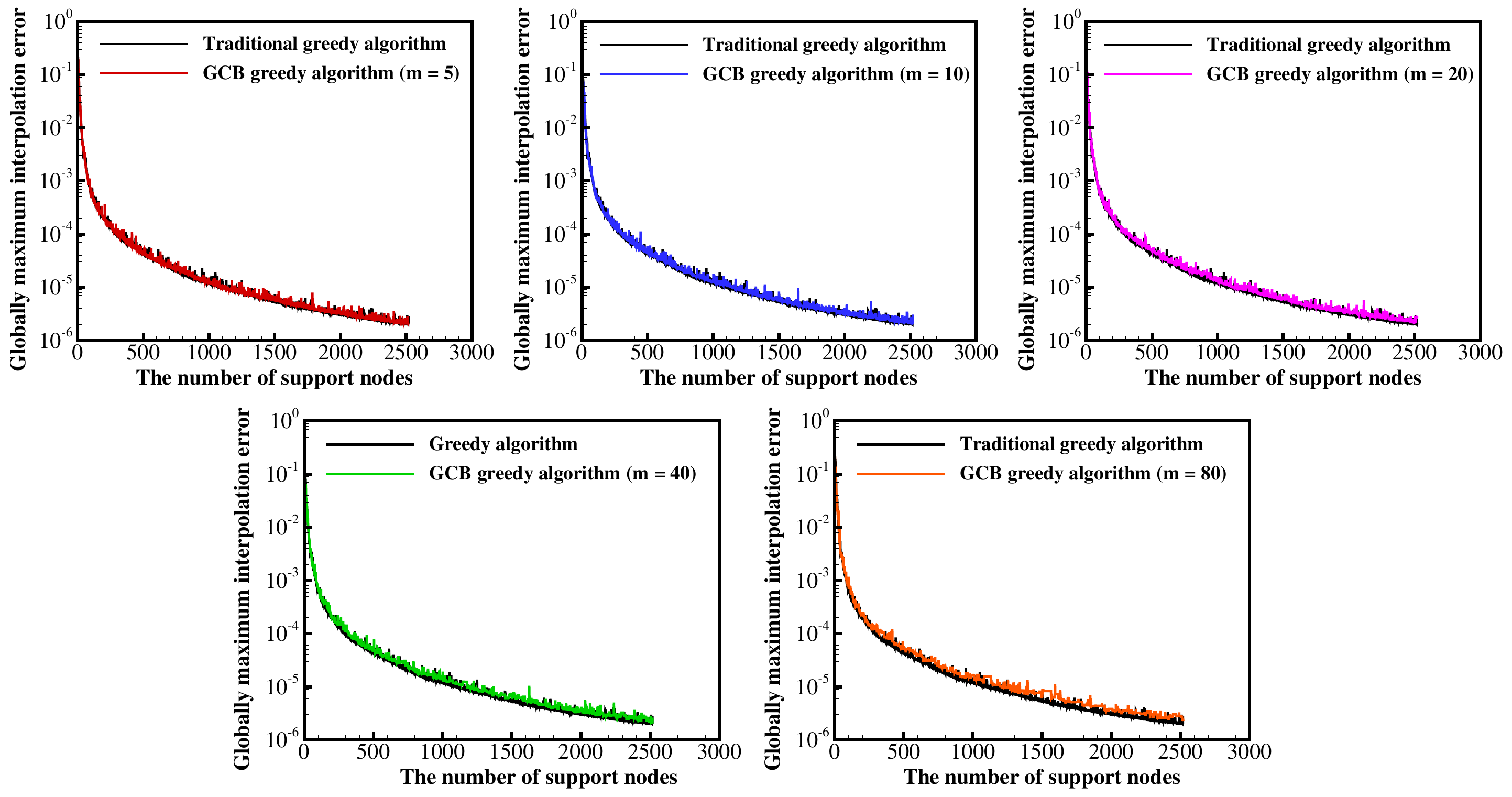}
\caption[1]{Globally maximum interpolation error histories in terms of the number of support nodes for the deformation of the DLR-F6 Wing-Body-Nacelle-Pylon configuration ($N_c = 2,520$).}\label{fig:13}
\end{figure}

\begin{figure}[htp]
\centering
\includegraphics[width=16cm]{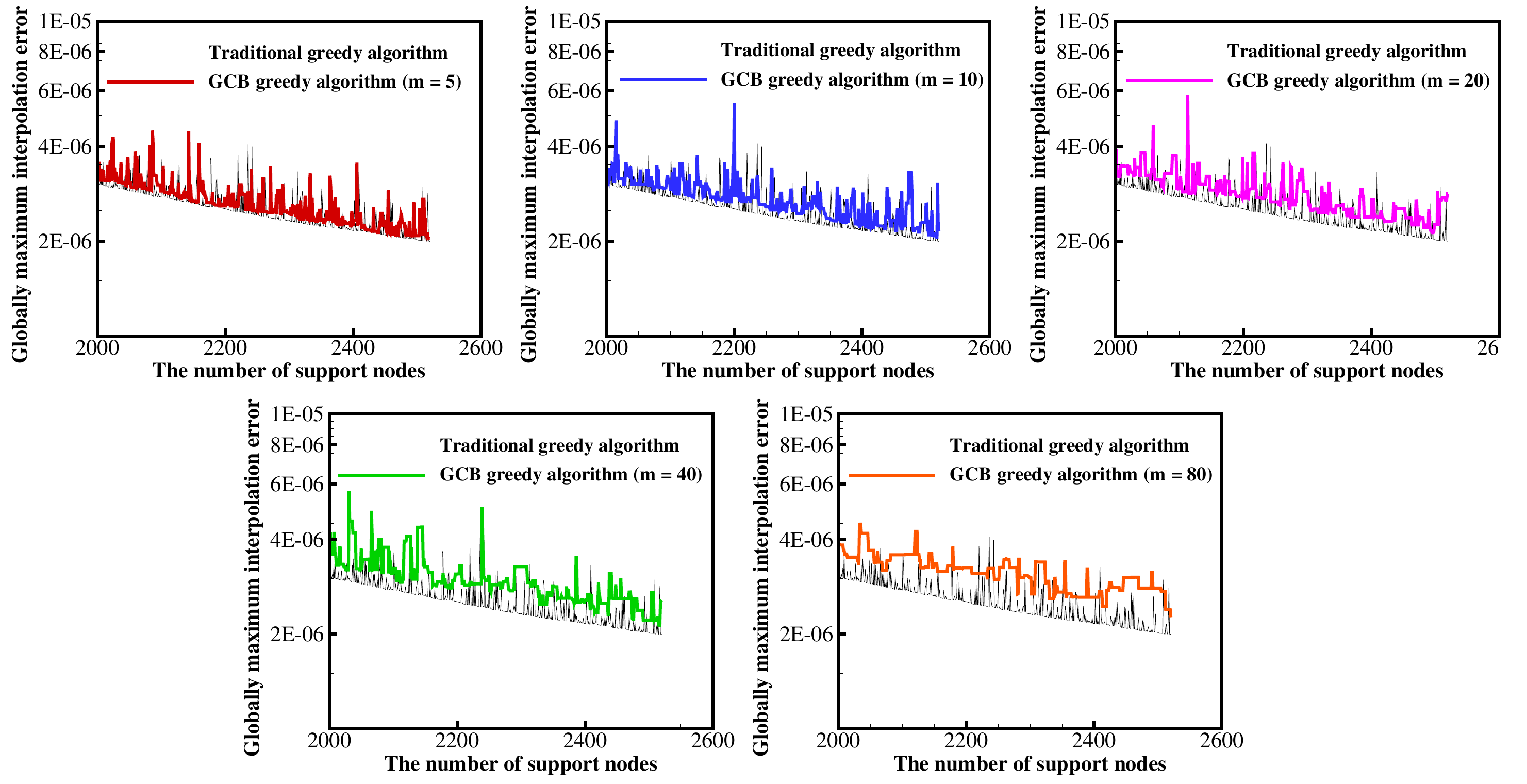}
\caption[1]{Globally maximum interpolation error histories in terms of the number of support nodes for the deformation of the DLR-F6 Wing-Body-Nacelle-Pylon configuration in a locally enlarged drawing ($N_c = 2,520$).}\label{fig:14}
\end{figure}

\begin{figure}[htp]
\centering
\includegraphics[width=9cm]{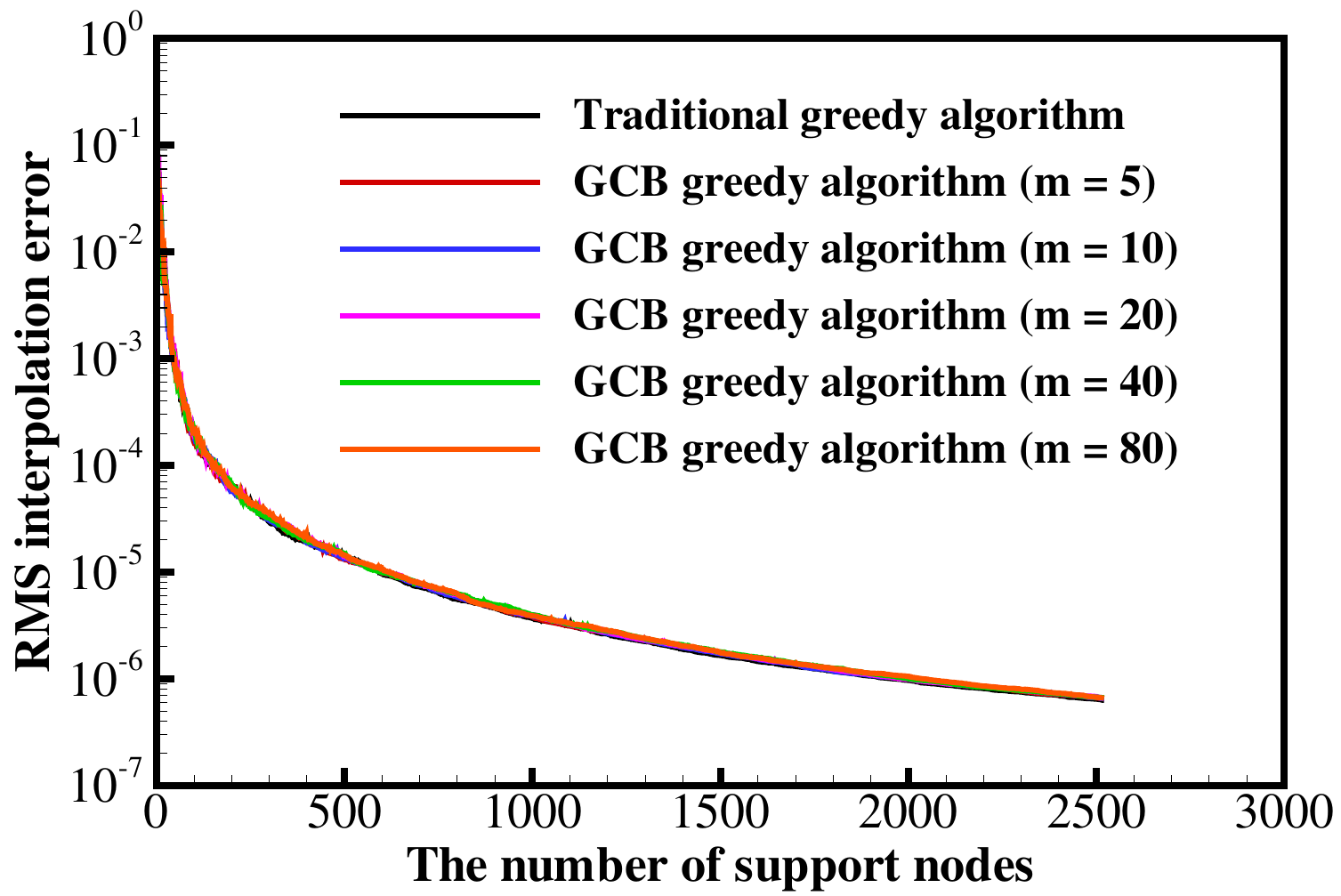}
\caption[1]{RMS interpolation error histories in terms of the number of support nodes for the deformation of the DLR-F6 Wing-Body-Nacelle-Pylon configuration ($N_c = 2,520$).}\label{fig:15}
\end{figure}

By allowing $N_c$ to vary, Table \ref{tab:6} shows the influences of $m$ on the time consumptions when performing mesh deformation. It is clearly seen that an increase of $m$ results in an increase of $N_c$ from 2520 to 2773 and a remarkable decrease of the time consumption for computing the interpolation errors ($t_1$) from \SI{357.24}{s} to \SI{6.49}{s}. However, it tends to generate increasing time consumptions for solving the linear algebraic system ($t_2$) and computing the displacements of volume nodes ($t_3$) in the case of $m > 50$, since the computational complexities for these two processes increases as functions of $N_c^3$ and $N_c$, respectively. For this reason, the efficiency of mesh deformation indicated by the sum of $t_1$, $t_2$ and $t_3$ doesn't keep increasing as $m$ grows. According to Table \ref{tab:6}, this sum decreases from \SI{393.6}{s} to the minimum of \SI{46.2}{s} when $m$ grows to 50, which is in the range of $\left[ {{N_b}/{N_c},{\rm{ }}2{N_b}/{N_c}} \right]$ ($N_b/N_c \approx 31.7$) suggested in section \ref{sec:4.1}, and the efficiency of mesh deformation is therefore improved by 8.5 times that is almost one order of magnitude. Compared to the case as shown in Table \ref{tab:3}, it is manifested that the use of the GCB greedy algorithm tends to generate a more significant efficiency improvement for mesh deformation when a larger-scale mesh is applied.

\begin{table}[ht!]\small
\caption{Number of support nodes and time consumptions for the deformation of the DLR-F6 Wing-Body-Nacelle-Pylon configuration.}\label{tab:6}
\begin{center}
\renewcommand{\arraystretch}{1.0}
\begin{tabular}{cccccc}
\hline \noalign{\smallskip}
Algorithm &$N_c$  & $t_1(s)$ & $t_2(s)$  & $t_3(s)$ & $t_1+t_2(s)+t_3(s)$ \\
\hline
Traditional greedy algorithm  &2520	&357.24   &5.51   &30.85   &393.60 \\
GCB greedy algorithm ($m = 5$)  &2544	&69.32    &5.57   &31.45   &106.48 \\
GCB greedy algorithm ($m = 10$) &2577	&34.78    &6.13   &32.03   &72.94 \\
GCB greedy algorithm ($m = 20$) &2567 &17.99    &5.92   &31.85   &55.76 \\
GCB greedy algorithm ($m = 40$) &2640 &10.10    &6.43   &32.43   &48.96 \\
GCB greedy algorithm ($m = 50$) &2655 &7.99     &5.58   &32.63   &46.20 \\
GCB greedy algorithm ($m = 80$) &2773 &6.49     &7.44   &34.10   &48.03 \\
\hline
\end{tabular}
\medskip
\end{center}
\end{table}

\begin{figure}[htp]
\centering
\subfigure{
\begin{minipage}[t]{0.52\linewidth}
\centering
\includegraphics[height=6.8cm]{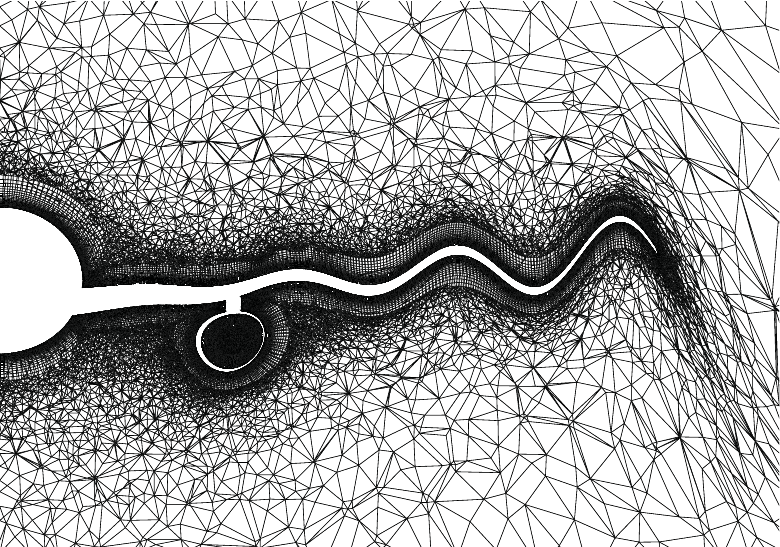}
\end{minipage}
}%
\subfigure{
\begin{minipage}[t]{0.46\linewidth}
\centering
\includegraphics[height=6.8cm]{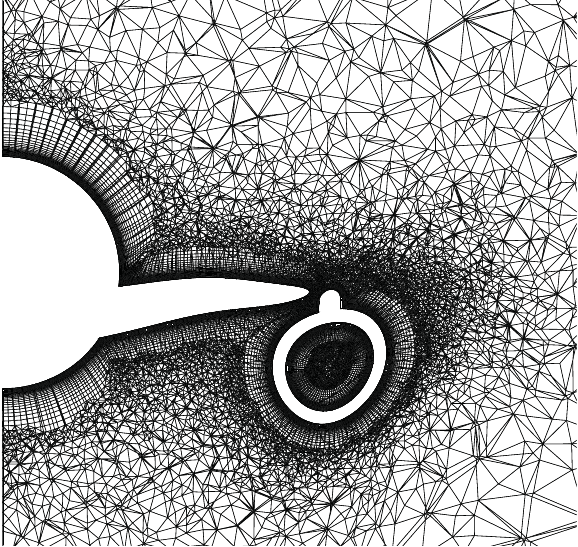}
\end{minipage}
}%
\centering
\caption{Mesh distribution of the DLR-F6 Wing-Body-Nacelle-Pylon configuration after deformation.}\label{fig:16}
\end{figure}

\begin{figure}[htp]
\centering
\includegraphics[width=12cm]{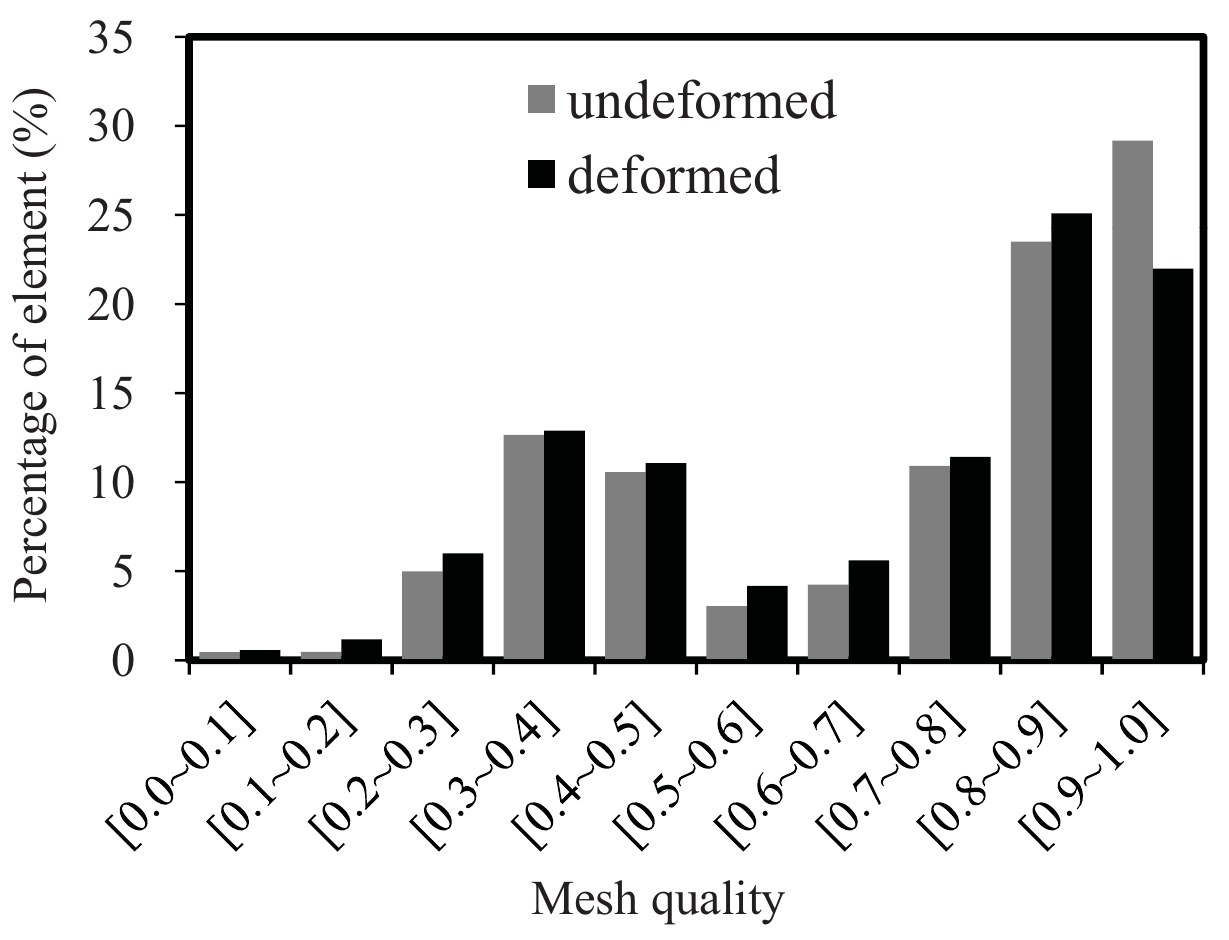}
\caption[1]{Mesh qualities of the DLR-F6 Wing-Body-Nacelle-Pylon configuration undeformed and deformed.}\label{fig:17}
\end{figure}

The mesh distribution after deformation is given in Fig. \ref{fig:16} and a comparison of the mesh qualities before and after deformation is given in Fig. \ref{fig:17}. It is shown that the GCB greedy algorithm can produce a deformed mesh with a comparable quality to the undeformed one for complex problems involving large deformations, which contributes to the subsequent aerodynamic simulation based on the computational fluid dynamics (CFD) technique.

\section{Conclusions}

In present work, a GCB greedy algorithm aimed at improving the efficiency of mesh deformation is proposed and validated in deformations of the ONERA M6 wing and the DLR-F6 Wing-Body-Nacelle-Pylon configuration. The conclusions are made in the following.
\begin{enumerate}

\item By incorporating the multigrid concept, the GCB greedy algorithm treats the maximum interpolation error of each group as an approximation to that of all boundary nodes and uses the node with the maximum interpolation error in each group to construct the set of support nodes. Unlike other boundary node-based reducing algorithms, it allows all boundary nodes can contribute to error control, thus ensuring the accuracy.

\item The computational results indicate that the GCB greedy algorithm is able to reduce the computational complexity for computing the interpolation errors in the data reducing procedure from $O\left( {N_c^2{N_b}} \right)$ to $O\left( {N_c^3} \right)$ and therefore promote the efficiency of this process by dozens of times. In addition, $m$ should be appropriately set in the range of $\left[ {{N_b}/{N_c},{\rm{\ }}2{N_b}/{N_c}} \right]$  to prevent too much additional computations for solving the linear algebraic system and computing the displacements of volume nodes induced by the increase of $N_c$ which is resulted from the increase of $m$.

\item It is manifested that the superiority of the GCB greedy algorithm in improving the efficiency of mesh deformation is more distinct for the use of a larger-scale mesh. For the deformation of the DLR-F6 Wing-Body-Nacelle-Pylon configuration, the efficiency of mesh deformation is improved by almost one order of magnitude.

\item For both structured and unstructured meshes, the GCB greedy algorithm can ensure the mesh quality after deformation and retain the grid orthogonality and the grid spacing near surface, which can therefore contribute to the subsequent aerodynamic simulation based on the CFD technique.
\end{enumerate}

\section*{Acknowledgement}
The financial assistances provided by the National Natural Science Foundation of China (No. 61773068) and the National Key Research and Development Program of China (No. 2019YFB1704204) are gratefully acknowledged.

\bibliography{mybibfile}
\end{document}